\documentclass{article}
\usepackage{amsmath,amssymb,pstricks,wasysym,stmaryrd,enumerate,epsfig}

\usepackage[french]{babel}
\newcommand{\tmem}[1]{{\em #1\/}}
\newcommand{\tmop}[1]{\operatorname{#1}}
\newcommand{\tmstrong}[1]{\textbf{#1}}
\catcode`\à=\active \defà{\`a} \catcode`\À=\active \defÀ{\`A}
\catcode`\á=\active \defá{\'a} \catcode`\Á=\active \defÁ{\'A}
\catcode`\ä=\active \defä{\"a} \catcode`\Ä=\active \defÄ{\"A}
\catcode`\â=\active \defâ{\^a} \catcode`\Â=\active \defÂ{\^A}
\catcode`\å=\active \defå{{\aa}} \catcode`\Å=\active \defÅ{{\AA}}
\catcode`\ç=\active \defç{\c{c}} \catcode`\Ç=\active \defÇ{\c{C}}
\catcode`\è=\active \defè{\`e} \catcode`\È=\active \defÈ{\`E}
\catcode`\é=\active \defé{\'e} \catcode`\É=\active \defÉ{\'E}
\catcode`\ë=\active \defë{\"e} \catcode`\Ë=\active \defË{\"E}
\catcode`\ê=\active \defê{\^e} \catcode`\Ê=\active \defÊ{\^E}
\catcode`\ì=\active \defì{\`{\i}} \catcode`\Ì=\active \defÌ{\`{\I}}
\catcode`\í=\active \defí{\'{\i}} \catcode`\Í=\active \defÍ{\'{\I}}
\catcode`\ï=\active \defï{\"{\i}} \catcode`\Ï=\active \defÏ{\"{\I}}
\catcode`\î=\active \defî{\^{\i}} \catcode`\Î=\active \defÎ{\^{\I}}
\catcode`\ò=\active \defò{\`o} \catcode`\Ò=\active \defÒ{\`O}
\catcode`\ó=\active \defó{\'o} \catcode`\Ó=\active \defÓ{\'O}
\catcode`\ö=\active \defö{\"o} \catcode`\Ö=\active \defÖ{\"O}
\catcode`\ô=\active \defô{\^o} \catcode`\Ô=\active \defÔ{\^O}
\catcode`\ù=\active \defù{\`u} \catcode`\Ù=\active \defÙ{\`U}
\catcode`\ú=\active \defú{\'u} \catcode`\Ú=\active \defÚ{\'U}
\catcode`\ü=\active \defü{\"u} \catcode`\Ü=\active \defÜ{\"U}
\catcode`\û=\active \defû{\^u} \catcode`\Û=\active \defÛ{\^U}
\catcode`\ý=\active \defý{\'y} \catcode`\Ý=\active \defÝ{\'Y}
\catcode`\ÿ=\active \defÿ{\"y} \catcode`\˜=\active \def˜{\"Y}
\catcode`\½=\active \def½{!`}
\catcode`\¾=\active \def¾{?`}
\catcode`\ß=\active \defß{{\ss}}

\begin{document}

\title{L'invariant de Bieri-Neumann-Strebel des groupes fondamentaux des
variétés kähleriennes}\author{Thomas Delzant}\maketitle

\section{Introduction.}

L'invariant de Bieri Neumann Strebel d'un groupe $G$de type fini, est le
sous-ensemble de $H^1 ( G, \mathbb{R} ) - \{ 0 \} / \mathbb{R}^{+ \ast}$ dont
le complémentaire, appelé ci-après $E^1 ( G ),$est constitué des classes
{\tmem{exceptionnelles}}, dont nous rappellerons la définition au moment
opportun.

Le but de cet article est de décrire cet invariant pour les {\tmem{groupes de
Kähler}}, c'est-à-dire les groupes fondamentaux des variétés kähleriennes
compactes. Les exemples les plus importants de variétés kähleriennes compactes
sont les variétés projectives lisses, et c'est une question importante de
comprendre la structure de leurs groupes fondamentaux (voir [ABC] pour une
introduction approfondie à ce sujet).

\subsection{Théorème.} {\tmem{Soit $X$ une variété kählerienne compacte et  $w
\in H^1 ( G ) .$ Alors $w \in E^1 ( G )$ si et seulement si il existe une
application holomorphe à fibres connexes $f : X \rightarrow S$, où $S$ est une
orbi-surface de Riemann hyperbolique de genre $g \geqslant 1$, et une 1-forme
fermée holomorphe sur $S$ $\omega$telle que $w = [ \tmop{Re} f^{\ast} \omega ]
.$}}

Rappelons qu'une orbi-surface de Riemann est une surface de Riemann de genre
$g$ marqués par $n$ points $p_1, \ldots . p_n$ équipés d'entiers $m_i
\geqslant 2$. Une telle surface est dite hyperbolique si sa caractéristique
d'Euler  $\chi^{\tmop{orb}} ( S ) = 2 - 2 g - \Sigma ( 1 - \frac{1}{m_i} )$
est strictement négative. C'est alors le quotient du disque unité par l'action
d'un réseau discret co-compact de $\tmop{PSL} ( 2, \mathbb{R} )$ ayant
$\mathbb{Z} / m_i \mathbb{Z}$ comme groupes d'isotropie aux points se
projetant sur les $p_i .$ Une application $X \rightarrow S$ est holomorphe si
elle se relève en un application holomorphe équivariante du revêtement
universel $\tilde{X}$ de $X$ vers le disque unité ; c'est le cas si et
seulement si elle l'est au sens usuel et si pour tout $i$ chaque composante
connexe de la fibre $f^{- 1} ( P_i )$ a une multiplicité  divisible par $m_i
.$ Les orbi-surfaces de Riemann de genre $g \geqslant 1$ ont toutes des
1-formes dont le noyau n'est pas de type fini et sont donc exceptionnelles au
sens de [BNS].

D'où l'idée d'utiliser l'invariant de Bieri Neumann Strebel (ou plutôt son
complémentaire $E^1 ( G ))$pour affiner ce que nous connaissons des groupes de
Kähler. Pour une introduction à la théorie de l'invariant BNS, voir  [L] ;
pour une étude approfondie [BS].  Rappelons en particulier le :

\subsection{Théorème [BNS],(voir aussi [BN] thm I.4.2)}{\tmem{ Si $\chi : G
\rightarrow A$ est un quotient abélien de $G$, $\ker \chi$ n'est pas de type
fini, si et seulement si il existe une classe {\tmem{exceptionnelle}} $w$
telle que $\ker w \supset \ker \chi .$}}

Il est convenu ([ABC]) de dire qu'une variété kählerienne $X${\tmem{fibre}} si
il existe une application holomorphe surjective à fibre connexe sur une
orbi-surface de Riemann hyperbolique. Il est bien connu que le fait qu'une
variété fibre ou non ne dépend que de son groupe fondamental. Dans le cas de
fibrations sur des surfaces de Riemann de genre $g \geqslant 2,$cela résulte,
par exemple, du théorème de Castelnuovo de Franchis (voir [ABC]) ; dans notre
langage une variété Kählerienne fibre sur une orbi-surface de Riemann
hyperbolique de genre $\geqslant 1$ si et seulement si $E^1 ( G ) \neq
\emptyset$.

En mélangeant les théorèmes  1.1 et 1.2, on répond à une question d'A.
Beauville [B] :

\subsection{Corollaire. } {\tmem{Si $G$ est le groupe fondamental d'une
variété kählerienne compacte $X$, A un groupe abélien et $\chi : G \rightarrow
A$ un homomorphisme dont le noyau n'est pas de type fini, alors $X$ fibre.
Plus précisément il existe une une orbi-surface de Riemann hyperbolique $S$,
une application holomorphe à fibres connexes $f : X \rightarrow S$,  et une
1-forme fermée holomorphe sur $S$ $\omega$telle que le noyau de la classe  $w
= [ \tmop{Re} f^{\ast} \omega ]$ contienne $\ker \chi .$}}

Notons que le cas particulier où $A$ est le groupe infini cyclique
$\mathbb{Z}$  et $w$ est un élément de $H^1 ( G,\mathbb{Z} )$ avait été
compris depuis longtemps par T. Napier et M. Ramachandran [NR].

Nous appliquerons ce corollaire à l'étude des groupes de Kähler, plus
particulièrement à l'étude de leurs quotients résolubles. En le combinant avec
les travaux de J.R. Groves [Gr] sur la structure des groupes résolubles,  nous
généraliserons au paragraphe 3 les résultats que D. Arapura, M. Nori [AN] et
F. Campana ([Ca]) avaient obtenu pour les groupes résolubles{\tmem{
linéaires.}}

\subsection{Théorème.} {\tmem{Si $G$ est le groupe fondamental d'une variété
kählerienne compacte $X$, et $p : G \rightarrow R$ un quotient résoluble de
$G$, non virtuellement nilpotent. Il existe revêtement fini $X_1 $de $X$, une
application holomorphe $f : X_1 \rightarrow S$ de $X_1$sur une orbi-surface de
Riemann hyperbolique et un quotient $R_1 $ résoluble non virtuellemnt
nilpotent de $R$ telle que l'application $\pi_1 ( X_1 ) \rightarrow R_1$ se
factorise par $f_{\ast}$. }} {\tmem{En particulier un groupe de Kähler
résoluble est virtuellement nilpotent.}}

\section{Démonstration du théorème 1.1}.

Partant de $w \in H^1 ( G, \mathbb{R} )$, la théorie de Hodge nous dit que $w$
peut être représentée par la partie réelle d'une 1-forme holomorphe $\omega$
de type $( 1, 0 )$ sur $X .$ Si $\pi : \tilde{X} \rightarrow X$ est un
revêtement universel de $X,$ $\pi^{\ast} \omega = d F$, où $F$ est une
fonction holomorphe, ce qui permet de définir les feuilles de $\omega$ comme
images par $\pi$ des composantes connexes des variétés de niveau de $F$. Il
s'agit de montrer que sous l'hypothèse que $w$ est exceptionnelle, les
feuilles de $w$ sont {\tmem{compactes}}, et sont les fibres d'une application
holomorphe sur une orbi-surface de Riemann.

Le lieu critique $C_w$ est le lieu des zéros de $w$. C'est un sous-ensemble
analytique complexe compact de $X$. Chaque composante connexe de ce
sous-ensemble est contenu dans une feuille du feuilletage. Appelons
{\tmem{feuille singulière}} une feuille de $w$ qui est une composante connexe
(de codimension complexe $1 )$ de $C_w$. L'existence d'une telle feuille
permet de conclure : c'est une feuille compacte. Or le fait suivant est connu
:

\subsection{Proposition. }([DG] 4.1 qui reprend un argument de [GS]).
{\tmem{Soit $X$ une variété Kählerienne compacte, $\omega$ une 1 forme
holomorphe fermée sur $X$. Si le feuilletage défini par $\omega$ a une feuille
compacte, alors il existe une orbi-surface de Riemann compacte  une
application holomorphe $X \rightarrow S$ dont les fibres sont les feuilles de
$\omega .$ }}

Nous devons donc nous convaincre que le fait que $w \in E^1$ entraîne
l'existence d'une feuille compacte {\tmem{singulière}}. Nous noterons toujours
$\pi : \tilde{X} \rightarrow X$ un revêtement universel, $F$ une primitive de
$\omega$ sur $\tilde{X}$, et  $h = \tmop{Re} F$ qui est une application
harmonique $G$ équivariante $h : \tilde{X} \rightarrow \mathbb{R}$.

\subsection{Lemme.}{\tmem{ Si}} $w \in E^1 ( G ),$ $h^{- 1} [ 0, + \infty [$
{\tmem{a au moins 2 composantes connexes non compactes.}}

Démonstration. On peut penser à $X$ comme à un CW complexe avec $X^0$ réduit à
un point et $X^1$ un bouquet de cercles. Le $1 -$squelette de $\tilde{X}$ est
le revêtement universel de $X$ et la non connexité de  l'intersection de $h^{-
1} [ 0, + \infty [$ avec ce 1-squelette est la {\tmem{définition}} de $E^1 ( G
) .$ $\Box$

On présente alors deux arguments pour prouver 1.1.

\subsubsection{Premier argument.} Grâce à 2.2., on observe que pour $t$
suffisamment grand $h^{- 1} ( t )$ ne peut être connexe. Le résultat principal
de C. Simpson [Si] s'applique.  $\Box$

\subsubsection{Second argument}.   Supposons que $C_w$ n'a pas de feuille
singulière et montrons, en raisonnant par l'absurde que $h^{- 1} [ 0, + \infty
[$ a une seule composante connexe non compacte. Considérons une valeur
régulière $r \gg 1$ de $h$, deux points $a, b$ qui ne sont pas dans la même
composante connexe de $h^{- 1} [ r, + \infty [$, et tels que $h ( a ), h ( b )
= r .$Le principe du minimax permet de construire une valeur critique $\xi$ :
pour tout chemin $\gamma$ joignant $a$ et $b$  dans $\tilde{X}$ considérons la
valeur minimale $m ( \gamma )$ de $h$ sur $\Gamma$. Soit $\xi = \max_{\gamma}
m ( \gamma )$, de sorte $\xi$ est  une valeur critique, et soit $\gamma : [ -
1, 1 ] \rightarrow ( M, a, b )$ un chemin critique avec $m ( \gamma ) = \xi =
h ( \gamma ( 0 ))$ et $h ( \gamma ( t )) > \xi$ si $t \neq 0 .$ Posons $c =
\gamma ( 0 )$, et montrons que la feuille de $\omega$ passant en $c$  est une
composante connexe de $\tilde{C}_w $, c'est à dire une feuille singulière.

Étudions d'abord un voisinage de $c.$ Dans une petite boule $B ( c, r )$,
identifiée à la boule unité de $\mathbb{C}^n$, on pense à $\tilde{C}_{\omega}$
comme à un germe de sous variété analytique complexe. Supposons par l'absurde
qu'elle est de codimension complexe $\geqslant 2$, réelle 4. En coupant par un
sous-espace vectoriel complexe générique de dimension 2 passant à l'origine,
on se ramène au cas où $\widetilde{C_{}}_{\omega}$ est de dimension 0
constituée de points singuliers isolés de la primitive $F$ de $\omega$.
Celle-ci est alors une fonction holomorphe ayant un point singulier isolé à
l'origine. Pour $\alpha$ suffisamment petit on peut alors considérer la sphère
$S_{\alpha}$ centrée en $c$ de rayon $\alpha$ et la fibration de Milnor de
$S_{\alpha} - ( C \cap S_{\alpha} ) \rightarrow S^1$. Or [Mi], celle-ci est à
fibres connexes car sa restriction au bord d'un petit voisinage tubulaire de 
de $C \cap S_{\alpha}$ est triviale (voir [W], chap.6 pour une démonstration
élémentaire). En particulier l'ensemble$( h - \xi ) > 0 \bigcap S_{\alpha} - (
C \cap S_{\alpha} )$ est {\tmem{connexe}}. Ainsi on peut joindre deux points
génériques de $h^{- 1} ( \xi + \varepsilon )$ proche de $c$ par un chemin
contenu dans $h ( x ) > \xi$, ce qui contredit le fait que $\xi$ est le
minimax. Donc au voisinage de $c$, $\tilde{C}_{\omega} \cap F^{- 1} ( F ( c
))$ est la feuille de $\omega$ passant en $c$. L'ensemble $\tilde{C}_{\omega}
\cap F^{- 1} ( F ( c ))$ est donc ouvert dans $F^{- 1} ( F ( c ))$, et il est
fermé car $\tilde{C}_{\omega}$ est le lieu critique de $F$ : $
\tilde{C}_{\omega} \cap F^{- 1} ( F ( c ))$ est la feuille de $\omega$ passant
en $c .\Box$

\section{Quotients résolubles des groupes de Kähler .}

Les résultats de 2 permettent d'affiner ce que nous connaissons sur la
structure des groupes de Kähler. Comme nous le notions déjà dans
l'introduction, le fait qu'une variété Kählerienne fibre ne dépend que de son
groupe fondamental.

\subsection{Lemme.}  {\tmem{Soit $G$ un groupe de Kähler. Supposons que $E^1
( G )^{} = \emptyset .$ Alors $G'$ est de type fini. De plus}}, {\tmem{le
métabélianisé $G / G^2$ de $G$ est virtuellement nilpotent.}}

{\tmstrong{Démonstration}}. Si  $E^1 ( G )^{} = \emptyset,$le groupe $G'$ est
de type fini. Donc son abélianisé $G'^{} / G^2$ aussi. Mais alors, d'après A.
Beauville [Be], il existe un nombre fini de caractères $\chi_1, \ldots \chi_k$
de $G$ dont les images sont des groupes finis cycliques, tels que les valeurs
propres de l'action de tout élément $g$ de $G$ sur $V = G'^{} / G^2 \otimes
\mathbb{C}$ soient de la forme $\chi_i ( g )$. L'action du sous groupe
d'indice fini $\cap \ker \chi_i$ de $G$ sur $V$ est donc unipotente.$\Box^{}$

Le théorème suivant est une conséquence des résultats de J.R.J. Groves [Gr],
récemment utilisés par  E. Breuillard [Br].

\subsection{Théorème.} {\tmem{ a) Soit $H$ un groupe résoluble de type fini
qui n'est pas virtuellement polycyclique. Alors $H$ admet un sous groupe
d'indice fini $H_1$ ayant un quotient $K_1$ métabélien  et tel que  $K_1' $ne
soit pas de type fini. }}

{\tmem{b) Soit $H$ est un groupe résoluble de type fini qui n'est pas
n'est pas virtuellement nilpotent. Alors $H$ admet un sous groupe d'indice
fini $H_1$ ayant un quotient $K_1$ métabélien qui n'est pas virtuellement
nilpotent.    }}

Démonstration. a) Considérons, avec J. Groves, un quotient $K$ de $H$ qui est
juste non virtuellement polycyclique, c'est-à-dire dont tout quotient propre
est virtuellement polycyclique. Rappelons qu'un tel quotient existe d'après le
lemme de Zorn car si $N_1 \subset N_2 \subset N_{k \ldots .}$ est une suite
croissante de sous-groupes distingués de $H$ tels que les quotients $H / N_i$
soient virtuellement polycycliques, la suite est stationnaire : en effet $H /
\cup N_i$ est virtuellement polycyclique de type fini donc de présentation
finie. Si $A \subset K$ est le sous-groupe de Fitting de $K$, Groves montre
que $A$ est abélien, que le quotient $K / A$ contient son centre $C$comme
sous-groupe d'indice fini, et $K$ est donc virtuellement métabélien. Le groupe
$A$ne peut être de type fini, sinon, $K$ serait virtuellement polycyclique.
L'image réciproque $H_1$ de $C$ est donc un groupe métabélien de type fini
extension d'un groupe abélien de type fini par un groupe abélien qui n'est pas
de type fini, et $H_1$ convient.

b) Considérons un quotient $K$ de $H$ qui est juste non virtuellement
nilpotent. Si $A$ est le sous-groupe de Fitting de $K$ Groves (voir aussi
l'article d'E. Breuillard [Br]) montre que $A$ est abélien et que $K / A$
contient son centre comme sous-groupe d'indice fini : ce groupe $K$ est donc
virtuellement métabélien non virtuellement nilpotent. $\Box$

On peut alors démontrer :

\subsection{Théorème.} {\tmem{Si $G$ est le groupe fondamental d'une variété
kählerienne $X$, et $p : G \rightarrow R$ un quotient résoluble de $G$, non
virtuellement nilpotent. Il existe revêtement fini $X_1 $de $X$ et une
application holomorphe $f : X_1 \rightarrow S$ de $X_1$sur une orbi-surface de
Riemann hyperbolique, et un quotient $R_1 $ résoluble non virtuellemnt
nilpotent de $R$ telle que l'application $\pi_1 ( X_1 ) \rightarrow R_1$ se
factorise par $f_{\ast} .$}}

{\tmem{Démonstration}}. D'après 3.2, $R$ admet un quotient $K$ virtuellement
métabélien mais pas virtuellement nilpotent. Soit $L \subset K$ un sous-groupe
d'indice fini qui est métabélien mais pas virtuellement nilpotent, et $G_1$ le
sous-groupe d'indice fini de $G$ qui est l'image réciproque de $L$.

Distinguons deux cas :

1. Si le groupe dérivé $L'$ est de type fini, $L$ est polycyclique ; mais
c'est un quotient de $G_1 / G_1^2$ qui est virtuellement nilpotent d'après
3.1, donc $L$ est aussi virtuellement nilpotent, contradiction.

2. Si $L'$ n'est pas de type fini,  le noyau de $\chi : G_1 \rightarrow
L^{\tmop{ab}} = L / L'$ se surjecte sur $L'$ qui est un groupe abélien qui
n'est pas de type fini. Le théorème 1.1 s'applique. $\Box$

\subsection{Corollaire.} {\tmem{Un groupe de Kähler résoluble est
virtuellement nilpotent.}}

En effet, un sous-groupe d'indice fini d'un groupe résoluble ne peut se
surjecter sur le groupe fondamental d'une (orbi)-surface hyperbolique, car
celui-ci contient un groupe libre. $\Box$

\subsubsection{Exemples.}  O. Kharlampovich [K] a construit un groupe
résoluble de présentation finie dont le problème du mot n'est pas résoluble ;
ce groupe ne peut pas être résiduellement fini, ni a fortiori, linéaire. D.
Robinson et R. Strebel [RS] ont construit un groupe résoluble de présentation
finie non linéaire dont le groupe dérivé n'est pas de type fini ; de tels
groupes ne sont pas des groupes fondamentaux de variétés kähleriennes.

{\large {\tmstrong{Remerciements.{\large }}}} Je tiens à remercier Gilbert
Levitt et Mihai Damian qui m'ont expliqué les bases de l'invariant de Bieri
Neumann Strebel. 

\section{Bibliographie.}

[ABC] Amorós, J.; Burger, M.; Corlette, K.; Kotschick, D.; Toledo, D.
Fundamental groups of compact Kähler manifolds. Mathematical Surveys and
Monographs, 44. American Mathematical Society, Providence, RI, 1996

[AN] Arapura, D. ; Nori, M.  Solvable fundamental groups of algebraic
varieties and Kähler manifolds.  Compositio Math. 116 (1999), no. 2, 173--188.

[B] Beauville, A. Annulation du $H^1$ pour les fibrés en droites plats.
{\tmem{Complex algebraic varieties }}(Bayreuth, 1990), 1--15, Lecture Notes in
Math., 1507

[BS] Bieri, Robert ; Strebel, Ralph. {\tmem{Geometric invariants for discrete
groups}}, cours polycopié.

[BNS] Bieri, Robert; Neumann, Walter D.; Strebel, Ralph A geometric invariant
of discrete groups. Invent. Math. 90 (1987), no. 3, 451--477.

[Br] Breuillard, E. On uniform exponential growth for solvable groups,
preprint 2006.

[C] Campana, F. Ensembles de Green-Lazarsfeld et quotients résolubles des
groupes de Kähler.  J. Algebraic Geom. 10 (2001), no. 4, 599--622.

[DG] Delzant,T. Gromov, M. Cuts in Kähler groups, in {\tmem{Infinite Groups:
Geometric, Combinatorial and Dynamical Aspects}},  Progress in Mathematics
Series 248,  Birkhauser,  2005.

[GS] Gromov M, Schoen, R. Harmonic maps into singular spaces and p-adic
superrigidity for lattices in groups of rank one. Inst. Hautes Études Sci.
Publ. Math. No. 76 (1992), 165--246.

[Gr] Groves, J. R. J. Soluble groups with every proper quotient polycyclic.
Illinois J. Math. 22 (1978), no. 1, 90--95.

[K] Kharlampovic, O.  A finitely presented solvable group with unsolvable word
problem.  Izv. Akad. Nauk SSSR Ser. Mat. 45 (1981), no. 4, 852--873, 928.

[L] Levitt, G. $\mathbb{R}$-trees and the Bieri-Neumann-Strebel invariant.
Publ. Mat. 38 (1994), no. 1, 195-200

[Mi] J. Milnor. {\tmem{Singular points of complex hypersurfaces}}, Annals of
Mathematics Studies 61 Priceton University Press, 1968

[NR] Napier T., Ramachandran M.,  Hyperbolic Kähler manifolds and proper
holomorphic mappings to Riemann surfaces. Geom. Funct. Anal. 11 (2001), no. 2,
382--406.

[RS] D. Robinson, R. Strebel. Some finitely presented soluble groups. J.
London Math. Soc. (2) 26 (1982) p. 435-440.

[S] Simpson, C. Lefschetz theorems for the integral leaves of a holomorphic
one-form. Compositio Math. 87 (1993), no. 1, 99--113.

[W] C.T.C. Wall {\tmem{Singular Points of Plane Curves}} London Math. Soc.
Student texts 63 Cambridge University Press.

\end{document}